\documentclass[11pt]{article}
\newcommand{\be}{\begin{equation}}
\newcommand{\ee}{\end{equation}}
\begin{document}
\title{\textbf {Algorithmic construction of Shimura - Taniyama - Weil parametrization of elliptic curves over the rationals}}
\author{ H. Gopalakrishna Gadiyar and R. Padma\\
Department of Mathematics\\School of Advanced Sciences\\ V.I.T. University, Vellore 632014 INDIA\\E-mail: gadiyar@vit.ac.in, rpadma@vit.ac.in}
\date{~~}
\maketitle

\begin{abstract}
In this note we give an algorithm to explicitly construct the modular parametrization of an elliptic curve over the rationals given the Weierstrass function $\wp (z)$.

\noindent MSC2010: 11B68, 11F03, 11S40, 14H52, 14L05, 33E05
\end{abstract}
\section{Introduction} We state a straight forward algorithm to construct the modular parametrization of elliptic curves over the rational numbers. The algorithm is implicit in the papers of Honda \cite{Honda}, Clarke \cite{Clarke}, Goldfeld \cite{Goldfeld} and Buchstaber and Bunkova \cite{Buchstaber:Bunkova} described below. It consists of the following steps: 

\section{The algorithm} 
\begin{description}
\item{Step 1.} Given the equation of the elliptic curve 
\be
E: y^2=4x^3-g_2x-g_3 \label{eq:elliptic}
\ee
over $Q$, we have the parametrization $x=\wp (z;g_2,g_3), y =\wp '(z;g_2,g_3)$.
\item{Step 2.} The formal exponential $f_E(T)$ of the formal group corresponding to this curve has been constructed explicitly by Buchstaber and Bunkova.
\be
f_E(T)=-2\frac{\wp(T;g_2,g_3)}{\wp '(T;g_2,g_3)}
\ee
\item{Step 3.} As $f_E(f_L(T))=f_L(f_E(T))=T$ the formal logarithm $f_L(T)$ is got from $f_E(T)$ by Lagrange inversion. $f_L(T)$ is related to the $L$-series of the elliptic curve. This is implicit in the work of Honda and Clarke. If
\be
f_L(T)=\sum_{n=1}^\infty \frac{a(n)}{n} T^n 
\ee
then
\be
L_E(s) = \sum_{n=1}^\infty \frac{a(n)}{n^s} \, .
\ee
and vice-versa.
\item{Step 4.} Goldfeld explicitly gives the modular parametrization as 
\be
x=\alpha (z)=\wp (F(z))=\wp \left (\sum_{n=1}^\infty \frac{a(n)}{n} e^{2\pi i nz} \right )
\ee
and 
\be
y=\beta (z)=\wp '(F(z))=\wp '\left (\sum_{n=1}^\infty \frac{a(n)}{n} e^{2\pi i nz} \right )\, .
\ee
\end{description}
\section{Formal group law} A commutative one dimensional formal group law over a ring $A$ is a formal power series (\cite{Hazewinkel}, \cite{Silverman})
\be
F(t_1, t_2)=t_1+t_2+\sum_{i,j\ge 1} a_{ij} t_1^i t_2^j
\ee
such that the following conditions hold:
$$F(t,0)=t, F(t_1, t_2)=F(t_2, t_1)
$$
and
$$
F(t_1, F(t_2,t_3))=F(F(t_1,t_2),t_3)
$$
The formal group $F$, the formal exponential $f_E$ and the formal logarithm $f_L$ are connected by the relations 
\be
F(t_1, t_2)=f_E(f_L(t_1)+f_L(t_2))
\ee
and
\be
f_E(f_L(T))=f_L(f_E(T))=T \, .
\ee

\section{Elliptic Curve and its formal group} 
Let $E$ be an elliptic curve $E: y^2=4x^3-g_2x-g_3$ over $Q$. Let ${\cal L}=\left\{m\omega _1+ n \omega _2: m,n \in Z\right\}$ be the lattice spanned by the two fundamental periods $\omega _1, \omega _2$ such that 
\be
g_2=g_2(L)=60 \sum_{l \in {\cal L}-\left\{0\right\}}\frac{1}{l^4}
\ee
and 
\be
g_3 = g_3(L)=140 \sum_{l \in {\cal L}-\left\{0\right\}}\frac{1}{l^6} .
\ee
We have the Weierstrass parametrization 
 $x=\wp (z), y=\wp '(z)$ where $\wp (z) $ the Weierstrass $\wp $ - function  
\be
\wp (z)=\frac{1}{z^2} + 2 \sum_{n \ge 1} G_{2n+2}(E, \omega) \frac{z^{2n}}{(2n)!}
\ee
Here 
\be
G_k=\frac{(-1)^k(k-1)!}{2} \sum_{l \in {\cal L}-\left\{0\right\}}\frac{1}{l^k}
\ee
for $k \ge 4$. Note that $G_k=0$ whenever $k$ is an odd integer. Also all the $G_k$'s can be got in terms of $g_2$ and $g_3$. 

The Bernoulli -Hurwitz numbers have been defined by Katz\cite{Katz} as follows. If $E$ is an elliptic curve given by (\ref{eq:elliptic}), then
\be
BH_k = 2kG_k {\rm ~for~} k \ge 4; =0 {\rm ~if~} k {\rm ~is~odd.}
\ee

Clarke \cite{Clarke} defines the universal Bernoulli numbers by the formula
\be
\frac{T}{f_E(T)}=\sum_{k \ge 0} \hat{B}_k \frac{T^k}{k!}. \label{eq:universal}
\ee
When $f_E(t)=e^T-1 $, we get the classical Bernoulli numbers $B_k$.

For the elliptic curve $E$ in the Weierstrass form, Buchstaber and Bunkova \cite{Buchstaber:Bunkova} have used the theory of solitons to make explicit constructions of formal group laws.  
\be
F(t_1, t_2)= t_1 +t_2 -b~m ~\frac{2g_2 +3g_3m}{4-g_2m^2 -g_3m^3}
\ee
is the formal group law of the elliptic curve $E$, where for $i=1,~2$, $(x_i,~y_i)\in E$
\be
t_i=-2\frac{x_i}{y_i}, s_i=-\frac{2}{y_i},
\ee
and 
\be
m=\frac{s_1-s_2}{t_1-t_2}, ~b=\frac{t_1s_2-t_2s_1}{t_1-t_2} \, ,
\ee
and the function 
\be
f_E(T)=-2\frac{\wp(T;g_2,g_3)}{\wp '(T;g_2,g_3)}
\ee
is the exponential of the formal group. 

Hence the $\hat{B}_k $ can be got in terms of the Bernoulli - Hurwitz numbers $BH_k$ in the case of elliptic curves.

\section{Formal logarithm and the $L$-series of an elliptic curve} The formal logarithm $f_L(T)$ is got from $f_E(T)$ by Lagrange inversion. If  
\be
f_L(T)=\sum_{n=1}^\infty \frac{a(n)}{n} T^n
\ee
denotes the formal logarithm of an elliptic curve $E$, then the $L$- series of the elliptic curve is given by
\be
L_E(s) = \sum_{n=1}^\infty \frac{a(n)}{n^s} \, .
\ee
and vice-versa \cite{Honda}.

If we take $f_E(T)=e^T-1 $ in (\ref{eq:universal}), we get the classical Bernoulli numbers $B_k$ and the corresponding formal logarithm is $f_L(T)=\log (1+T)$ In this case one gets the Dirichlet ($L$-) series
\be
\eta (s) = 1-\frac{1}{2^s}+\frac{1}{3^s}-\cdots  = (1-2^{1-s})~\zeta (s) \, ,
\ee
where $\zeta (s)$ is the classical Riemann - zeta function.

The Shimura - Taniyama - Weil conjecture which is now known as the modularity theorem states that any elliptic curve over rational numbers is modular. 

In \cite{Goldfeld}, D. Goldfeld gives the modular parametrization of $E$ as follows:
\be
\alpha (z)=\wp (F(z))=\wp \left (\sum_{n=1}^\infty \frac{a(n)}{n} e^{2\pi i nz} \right )
\ee
and 
\be
\beta (z)=\wp '(F(z))=\wp '\left (\sum_{n=1}^\infty \frac{a(n)}{n} e^{2\pi i nz} \right )\, .
\ee
where
\be
f(z)=\sum_{n=1}^\infty a(n) e^{2\pi inz}
\ee
is a cusp form of weight 2 for $\Gamma_0(N)$ and $N$ is the conductor of the elliptic curve $E$.
\section{Conclusion} The Riemann zeta-function $\zeta (s)$ has two fundamental representations. 
\begin{equation}
\sum_{n=1}^\infty \frac{1}{n^s} ~=~ \prod_p \left(1-\frac{1}{p^s}\right)^{-1}, {\rm ~if~}Re. ~s~>~1 \,.
\end{equation}
These we call the left hand and the right hand side. The left hand side is a sum over the integers and the right hand side is a product over primes. In recent times the $p$-adic approach which corresponds to the right hand side has been cultivated extensively. It is hoped that the present note will re-kindle an interest in an approach which relies on complex analysis.

\section{Acknowledgements}  We would like to thank Professors Arul Lakshminarayanan (IIT, Madras), A. Sankaranarayanan (TIFR), K. Srinivas (IMSc, Chennai) for making available the papers and books used in this research.

\end{document}